\documentstyle[12pt]{article}

\newcommand{\rom}{\mathrm}

\newcommand{\N}{\rm{I\!N}}         
\newcommand{\R}{\rm{I\!R}}           
\newcommand{\Q}{\rm{I\!Q}}

\newtheorem{theorem}{Theorem}%[section]
\newtheorem{definition}[theorem]{Definition}

\newtheorem{coro}[theorem]{Corollary}

\newcommand{\eps}{\varepsilon}

\newcommand{\ohne}{\setminus}

\newcommand{\gen}{\rightarrow}

\newcommand{\aus}{\subset}
\newcommand{\Norm}[1]{\Bigl\|#1\Bigr\|}
\newcommand{\norm}[1]{\|#1\|}
\newcommand{\Betr}[1]{\Bigl| #1  \Bigr|}
\newcommand{\betr}[1]{| #1  |}
\newcommand{\eing}[1]{_{|{#1}}}
\newcommand{\ebew}{\hfill{\rule{1.2ex}{1.2ex}}}
\newcommand{\bgl}{\begin{eqnarray}}
\newcommand{\bglst}{\begin{eqnarray*}}
\newcommand{\egl}{\end{eqnarray}}
\newcommand{\eglst}{\end{eqnarray*}}
\newcommand{\Pel}{Pe\l\-czy\'ns\-ki}
\newcommand{\Ref}[1]{(\ref{#1})}

\newcommand{\eins}{{1}}
\newcommand{\id}{{{\rom i}\rom{d}}}
\newcommand{\charF}[1]{\eins_{#1}}
\newcommand{\wst}{$w^{*}$}

\begin{document}
\begin{center}
{\bf\Large A separable L-embedded Banach space has property (X)
and is therefore the unique predual of its dual.}\bigskip\\
{\large H.\ Pfitzner}
\end{center}
\begin{abstract}\noindent
In this note the following is proved.
Separable L-embedded spaces
- that is separable Banach spaces which are complemented in their biduals such
that the norm between the two complementary subspaces is additive - 
have property (X) which, by a result of Godefroy and Talagrand, entails
uniqueness of the space as a predual.
\end{abstract}
\vspace{1cm}
We say that a Banach space $X$ is the unique predual of its dual
(more precisely the unique isometric predual of its dual) in case
it is isometric to any Banach space whose dual is isometric to
the dual of $X$.
(We say that two Banach spaces $Y$ and $Z$ are isomorphic
if there is
a bounded linear bijective operator $T:Y\rightarrow Z$ with bounded
inverse $T^{-1}$; if moreover $\norm{T(y)}=\norm{y}$ for all $y\in Y$
we say that $Y$ and $Z$ are isometric.)
In general a Banach space need not be the unique predual of its dual,
for example $c$ and $c_0$ are not isometric Banach spaces although
their duals are.

As shown by Grothendieck \cite[Rem.\ 4]{Gro1} in 1955, $L^1$-spaces are unique
preduals of their duals. Using essentially a result of Dixmier \cite{Dix-Fo} from
1953, Sakai \cite[Cor.\ 1.13.3]{Sakai} observed that more generally preduals of von
Neumann algebras are unique, and Barton and Timoney \cite{Bar-Tim}
 and Horn \cite{Horn1} generalized this to preduals of $JBW^*$-triples.
Ando \cite{Ando-Hu} stated the uniqueness as a predual for the quotient
$L^1/H_0^1$.

As Banach spaces these examples have in common to be L-summands
in their biduals or, for short, to be L-embedded.
By definition a Banach space $X$ is L-embedded if there is a projection $P$ on
its bidual $X^{**}$ with range $X$ such that
$\norm{Px^{**}}+\norm{x^{**}-Px^{**}}=\norm{x^{**}}$ for all
$x^{**}\in X^{**}$.

The standard reference for L-embedded spaces is \cite{HWW},
for a survey on unique preduals we refer to \cite{G-survey},
for general Banach space theory to  \cite{JohLin}, \cite{LiTz12}, or
\cite{Die-Seq}.

If not stated otherwise a sequence $(z_j)$ (and similarly a series $\sum z_j$) is
indexed by $\N=\{1, 2, \ldots\}$; we write $\N_0=\{0, 1, 2, \ldots\}$.
Recall that a series $\sum z_j$ in a Banach space $Z$ is called
weakly unconditionally Cauchy (wuC for short) if
$\sum \betr{z^*(z_j)}$ converges for each $z^*\in Z^*$ or, equivalently,
if there is a number $M$ such that
$\norm{\sum_{j=1}^{n}\alpha_j z_j}\leq M\max_{1\leq j\leq n}\betr{\alpha_j}$
for all $n\in\N$ and all scalars $\alpha_j$.
It is well known by a result of Bessaga and \Pel\ that a Banach space contains
a subspace isomorphic to $c_0$ if and only if it contains a wuC-series
$\sum z_j$ such that $\inf\norm{z_j}>0$.
In case $\sum x_j^*$  is a wuC-series in a dual Banach space $X^*$
we denote the \wst-limit (that is the limit in the $\sigma(X^*,X)$-topology)
of the sequence $(\sum_{j=1}^{n} x_j^*)$  by $\sum^* x_j^*$.

In their study of unique preduals Godefroy and Talagrand \cite{GoTa}
defined \begin{definition}[Property (X)]
A Banach space $X$ has property (X) if for each $x^{**}\in X^{**}\ohne X$
there exists a wuC-series $\sum x_j^*$ in $X^*$ such that
\bgl
\sum x^{**}(x_j^*) \neq x^{**}\Bigl(\sum\mbox{}\!^* x_j^*\Bigr).\label{gl01}
\egl
\end{definition}
They proved
\begin{theorem}[Godefroy, Talagrand]\label{theo-GT}
A Banach space $X$ with property (X) is the unique predual of its dual.

Moreover, every $x^{**}\in X^{**}$ which is strongly Baire mesurable
on $(X^*, w^*)$ (for the definition see \cite[Th.\ V.3]{G-survey})
- in particular, every $x^{**}\in X^{**}$ which is Borel on $(X^*, w^*)$
- belongs to $X$.
\end{theorem}

Up to now it has been known that a class of L-embedded spaces,
namely the duals of M-embedded spaces (see \cite{HWW} for the definiton),
have property (X) \cite[p.\ 148]{HWW}. Furthermore it has been known (\cite{Pf-L1},
\cite[Th.\ VI.2.7]{HWW})
that L-embedded spaces have \Pel's property (V$^*$); the latter one is
similar to and implied by property (X) but is, by an example of
Talagrand \cite{Ta-Bor}, strictly weaker than (X).

In view of all this it was natural to ask whether L-embedded spaces are
unique preduals and have property (X)
(see \cite[Problem page 123]{HWW}).
At least in the separable case the answer is yes.
\begin{theorem}\label{theo3}
Separable L-embedded Banach spaces have property (X).
\end{theorem}
{\it Proof}:
Let $X$ be L-embedded and $P$ be the corresponding projection
on $X^{**}$ with range $X$; we put $Q=\id_{X^{**}}-P$.
Denoting the range of $Q$ by $X_{{\rm s}}$ we have the decomposition
$X^{**}=X\oplus_1 X_{{\rm s}}$.
Let the sequence $(x_n)$ be dense in $X$.
Let $x^{**}\in X^{**}\ohne X$.
Let $\eta=\norm{x_{{\rm s}}}$ where $x^{**}=x+x_{{\rm s}}$, $x\in X$,
$x_{{\rm s}}\in X_{{\rm s}}$.
We have that $\eta>0$ because $x^{**}\not\in X$.

Let $1>\eps>0$.
By the Bishop-Phelps theorem \cite{Bish-Ph}, \cite{Hol-FA} there
are $x^* \in X^*$ and $\tilde{x}^{**}\in X^{**}$ such that $\norm{x^*}=1$
and $\norm{\tilde{x}^{**} -x^{**}}<\eps\eta/3$ and such that $\tilde{x}^{**}$ attains its
norm on $x^*$ that is $\tilde{x}^{**}(x^*)=\norm{\tilde{x}^{**}}$.
For the decomposition $\tilde{x}^{**}=\tilde{x}+\tilde{x}_{{\rm s}}$ we have
\bglst
\tilde{x}_{{\rm s}}(x^*)=\norm{\tilde{x}_{{\rm s}}} \quad\mbox{ and }\quad
\norm{\tilde{x}_{{\rm s}}-x_{{\rm s}}}<\eps\eta/3
\eglst
because
$\norm{\tilde{x}_{{\rm s}}}
\geq
\betr{\tilde{x}_{{\rm s}}(x^*)}
=\betr{\tilde{x}^{**}(x^*)-x^*(\tilde{x})}
\geq\norm{\tilde{x}^{**}}-\norm{\tilde{x}}=\norm{\tilde{x}_{{\rm s}}}$
and
$\eps\eta/3 >\norm{\tilde{x}^{**}-x^{**}}=\norm{\tilde{x}-x+\tilde{x}_{{\rm s}}-x_{{\rm s}}}=
\norm{\tilde{x}-x}+ \norm{\tilde{x}_{{\rm s}}-x_{{\rm s}}}$.\medskip\\
Choose a sequence $(\eps_j)$ of strictly positive numbers
such that $\prod_{j=1}^{\infty} (1+\eps_j)<1+\eps$ and
$\prod_{j=1}^{\infty} (1-\eps_j)>1-\eps$.

By induction over $\N_0$ we construct two sequences
$(x_n^*)_{n\in\N_0}$ and $(y_n^*)_{n\in\N_0}$  in $X^*$
(of which the first members $x_0^*$ and $y_0^*$ are auxiliary  values used
only for the induction) such that,
for all (real or complex) scalars $\alpha_j$, the following holds:
\bgl
x_0^*=0, && \norm{y_0^*}=1,                                              \label{gl0}\\
y_n^*&=&  x^*-\sum_{j=0}^{n}x_j^*                                                      \label{gl1}\\
\Bigl(\prod_{j=1}^{n} (1-\eps_j)\Bigr)\max_{0\leq j\leq n}\betr{\alpha_j}
&\leq&
\Norm{\alpha_0 y_n^*+\sum_{j=1}^{n}\alpha_j x_j^*}
\nonumber\\
&\leq&
\Bigl(\prod_{j=1}^{n} (1+\eps_j)\Bigr)\max_{0\leq j\leq n}\betr{\alpha_j},
                                                      \quad\mbox{if } n\geq 1,       \label{gl2}\\
\tilde{x}_{{\rm s}}(x_j^*)&=&0 \quad\quad\quad\mbox{if }0\leq j\leq n,\label{gl3}\\
x_{{\rm s}}(x_j^*)&=&0 \quad\quad\quad\mbox{if }0\leq j\leq n,\label{gl4}\\
y_n^*(x_k)&=& 0 \quad\quad\quad\mbox{if }1\leq k\leq n.\label{gl5}
\egl
For $n=0$ we set $x_0^*=0$ and $y_0^*=x^*$.\\
We notice that the restriction of $P^*$ to $X^*$ is an isometric isomorphism
from $X^*$ onto $X_{{\rm s}}^{\bot}$,
that $Q$ is a contractive projection and that
$X^{***}=X_{{\rm s}}^{\bot}\oplus_{\infty} X^{\bot}$.\\
For the induction step $n\mapsto n+1$ suppose now that
$x_0^*, \ldots, x_n^*$ and $y_0^*, \ldots, y_n^*$ have been constructed. Put
\bglst
E&=&{\rm lin}(\{x^*, x_0^*, \ldots, x_n^*, y_n^*,  P^*x_0^*, \ldots, P^*x_n^*, 
P^*y_n^*\})\aus X^{***},\\
F&=&{\rm lin}(\{ x_1, \ldots, x_{n+1}, x_{{\rm s}}, \tilde{x}_{{\rm s}}\})\aus X^{**}.
\eglst 
Note that $Q^*x_j^*,\, Q^*y_n^*\in E$ for $0\leq j\leq n$.
By the principle of local reflexivity there is an operator
$R:E\rightarrow X^*$ such that
\bgl
(1-\eps_{n+1})\norm{e^{***}}&\leq&\norm{Re^{***}}
\leq(1+\eps_{n+1})\norm{e^{***}},\label{gl7}\\
f^{**}(Re^{***})&=&e^{***}(f^{**}),\label{gl8}\\
R\eing{E\cap X^*}&=&\id_{E\cap X^*}\label{gl9}
\egl
for all $e^{***}\in E$ and $f^{**}\in F$.

We define
\bglst
x_{n+1}^*=RP^*y_n^*  \,\,\,\,\mbox{ and }\,\,\,\,  y_{n+1}^*=RQ^*y_n^*.
\eglst
First we notice that (\ref{gl1}, $n+1$) holds because
\bglst
x^*-\sum_{j=0}^{n+1}x_j^* \stackrel{\Ref{gl1}}{=}y_n^*-x_{n+1}^* =R(y_n^*-P^*y_n^*) = RQ^*y_n^* =y_{n+1}^*.
\eglst
In the following we use the convention $\sum_{j=1}^{0} (\cdots) =0$. Then we have that
\bglst
\alpha_0 y_{n+1}^* +\sum_{j=1}^{n+1}\alpha_j x_j^*
= R\Bigl(Q^*(\alpha_0 y_n^* +\sum_{j=1}^{n}\alpha_j x_j^*) + P^*(\alpha_{n+1} 
y_n^* +\sum_{j=1}^{n}\alpha_j x_j^*)\Bigr).
\eglst
The second inequality of (\ref{gl2}, $n+1$) can be seen as follows:
\bgl
\Norm{\alpha_0 y_{n+1}^* +\sum_{j=1}^{n+1}\alpha_j x_j^*}
&\stackrel{\Ref{gl7}}{\leq}&
(1+\eps_{n+1})
\Norm{Q^*(\alpha_0 y_n^* +\sum_{j=1}^{n}\alpha_j x_j^*) +P^*(\alpha_{n+1} y_n^* 
+\sum_{j=1}^{n}\alpha_j x_j^*)}\nonumber\\
&=&
(1+\eps_{n+1})\max\Bigl\{\Norm{Q^*(\alpha_0 y_n^* +\sum_{j=1}^{n}\alpha_j 
x_j^*)},
                    \Norm{P^*(\alpha_{n+1} y_n^* +\sum_{j=1}^{n}\alpha_j 
x_j^*)}\Bigr\}\nonumber\\
&\leq&
(1+\eps_{n+1})\max\Bigl\{\Norm{\alpha_0 y_n^* +\sum_{j=1}^{n}\alpha_j x_j^*},
                    \Norm{\alpha_{n+1} y_n^* +\sum_{j=1}^{n}\alpha_j x_j^*}\Bigr\}\nonumber\\
&\leq&
\Bigl(\prod_{j=1}^{n+1}(1+\eps_j)\Bigr)\max\{\max_{0\leq j\leq n}\betr{\alpha_j},
                                             \max_{1\leq j\leq n+1}\betr{\alpha_j}\}\nonumber\\
&=&
\Bigl(\prod_{j=1}^{n+1}(1+\eps_j)\Bigr)\max_{0\leq j\leq n+1}\betr{\alpha_j}   \label{gl6}
\egl
where the last inequality comes from \Ref{gl0} if $n=0$ and from \Ref{gl2} if $n\geq 1$.\\
For the first inequality of (\ref{gl2}, $n+1$) we estimate
\bglst
\Norm{\alpha_0 y_{n+1}^* +\sum_{j=1}^{n+1}\alpha_j x_j^*}
&\stackrel{(\ref{gl7})}{\geq}&
(1-\eps_{n+1})\Norm{Q^*(\alpha_0 y_n^* +\sum_{j=1}^{n}\alpha_j x_j^*) +
                                  P^*(\alpha_{n+1} y_n^* +\sum_{j=1}^{n}\alpha_j x_j^*)}\\
&=&
(1-\eps_{n+1})\max\Bigl\{\Norm{Q^*(\alpha_0 y_n^* +\sum_{j=1}^{n}\alpha_j 
x_j^*)},
                    \Norm{\alpha_{n+1} y_n^* +\sum_{j=1}^{n}\alpha_j 
x_j^*}\Bigr\};
\eglst
in case $\betr{\alpha_0} =\max_{0\leq j\leq n+1}\betr{\alpha_j}$
we observe that $Q\tilde{x}_{{\rm s}}=\tilde{x}_{{\rm s}}$,
that $\tilde{x}_{{\rm s}}(y_n^*)=\tilde{x}_{{\rm s}}(x^*)-\sum_{j=1}^{n}\tilde{x}_{{\rm s}}(x_j^*)=
 \tilde{x}_{{\rm s}}(x^*)$ by \Ref{gl3},
and we continue the estimate by
\bgl
\cdots
&\geq&
(1-\eps_{n+1})\, \Betr{\Bigl(Q^*(\alpha_0 y_n^* 
      +\sum_{j=1}^{n}\alpha_j x_j^*)\Bigr)(\frac{\tilde{x}_{{\rm s}}}{\norm{\tilde{x}_{{\rm s}}}})}\nonumber\\
&=&
\frac{(1-\eps_{n+1})}{\norm{\tilde{x}_{{\rm s}}}}
\Betr{\tilde{x}_{{\rm s}} (\alpha_0 y_n^* 
      +\sum_{j=1}^{n}\alpha_j x_j^*)}\nonumber\\
&\stackrel{\Ref{gl3}}{=}&
\frac{(1-\eps_{n+1})}{\norm{\tilde{x}_{{\rm s}}}}
\, \betr{\alpha_0}\,
\betr{\tilde{x}_{{\rm s}}(y_n^*)}
=\frac{(1-\eps_{n+1})}{\norm{\tilde{x}_{{\rm s}}}}
\, \betr{\alpha_0}\,
\betr{\tilde{x}_{{\rm s}}(x^*)}\nonumber\\
&=&
(1-\eps_{n+1})\, \betr{\alpha_0}             \label{gl6bis}
\egl
whereas in case
$\betr{\alpha_0}\neq\max_{0\leq j\leq n+1}\betr{\alpha_j}$
we get
\bgl
\cdots
&\geq&
(1-\eps_{n+1})\, \Norm{\alpha_{n+1} y_n^* +\sum_{j=1}^{n}\alpha_j x_j^*}\nonumber\\
&\geq&
\Bigl(\prod_{j=1}^{n+1}(1-\eps_j)\Bigr)
\max_{1\leq j\leq n+1}\betr{\alpha_j}           \label{gl6bisbis}
\egl
where the last inequality comes from \Ref{gl0} if $n=0$ and from \Ref{gl2} if $n\geq 1$.
Thus we obtain the first inequality of (\ref{gl2}, $n+1$).

The conditions (\ref{gl3}, $n+1$), (\ref{gl4}, $n+1$) and (\ref{gl5}, $n+1$) are easy to
verify because $Px_{{\rm s}}=P\tilde{x}_{{\rm s}}=Qx_k=0$ thus
\bglst
x_{{\rm s}}(x_{n+1}^*)&=&x_{{\rm s}}(RP^*y_n^*)=(P^*y_n^*)(x_{{\rm s}})=Px_{{\rm s}}(y_n^*)=0,\\
\tilde{x}_{{\rm s}}(x_{n+1}^*)&=&\tilde{x}_{{\rm s}}(RP^*y_n^*)=(P^*y_n^*)(\tilde{x}_{{\rm s}})=P\tilde{x}_{{\rm s}}(y_n^*)=0,\\
y_{n+1}^*(x_k)&=&(RQ^*y_n^*)(x_k)=y_n^*(Qx_k)=0.
\eglst
This ends the induction.

By \Ref{gl2}, $\sum x_j^*$ is wuC (where, as indicated above in the
introduction, $j$ runs through $\N$).
We have that $\sum^* x_j^*=x^*$ by \Ref{gl1} (and \Ref{gl0})
because by \Ref{gl5} and the density of the $x_k$
the \wst-limit of $(y_n^*)$ is $0$.
This easily entails \Ref{gl01} because we have
$\sum x_{{\rm s}}(x_j^*)=0$ by \Ref{gl4}, we have
$\norm{\tilde{x}_{{\rm s}}}
\geq\norm{x_{{\rm s}}}-\norm{x_{{\rm s}}-\tilde{x}_{{\rm s}}}
>(1-\eps/3)\eta$
and trivially $(\sum^* x_j^*)(x) = \sum x_j^*(x)$ thus
\bglst
x^{**}(\sum\mbox{}\!^* x_j^*)-\sum x^{**}(x_j^*)
&=&
x_{{\rm s}}(\sum\mbox{}\!^* x_j^*)-\sum x_{{\rm s}}(x_j^*)\\
&=&
x_{{\rm s}}(x^*)=
\norm{\tilde{x}_{{\rm s}}}-(\tilde{x}_{{\rm s}}-x_{{\rm s}})(x^*)\\
&>&
(1-\frac{\eps}{3})\eta-\frac{\eps\eta}{3}\\
&>&
(1-\eps)\norm{Qx^{**}}>0.
\eglst
This ends the proof.\ebew\\
We have already mentioned Godefroy's and Talagrand's result
that property (X) implies the uniqueness of a Banach space as a predual;
moreover, since (X) is hereditary and stable by
equivalent norms we obtain
\begin{coro}
A Banach space that is isomorphic to a subspace of a separable
L-embedded space is the unique predual of its dual.
\end{coro}
Remarks:\\
1.
It follows immediately from the first variant of the proof of theorem \ref{theo3} that
if $\tilde{x}_{{\rm s}}$ is a non zero norm attaining element of
$X_{{\rm s}}$
then the two expressions in \Ref{gl01} differ by the greatest possible
value, more precisely
\bglst
\tilde{x}_{{\rm s}}(\sum\mbox{}\!^* x_j^*)
=\norm{\tilde{x}_{{\rm s}}}
\neq
0=\sum \tilde{x}_{{\rm s}}(x_j^*)
\eglst
with $\norm{\sum^* x_j^*}=1$.\\
2.
Is it possible to refine the proof of theorem \ref{theo3}
so to produce a sequence spanning $c_0$ almost or asymptotically
isometrically?
We say that a sequence $(z_j)$ in a Banach space $Z$ spans $c_0$ almost
isometrically if  there exists a sequence
$(\delta_m)$  satisfying $[0,1[\ni\delta_m\gen 0$ such that
$(1-\delta_m)\sup_{m\leq j\leq n}\betr{\alpha_j}
\leq
\Norm{\sum_{j=m}^{n} \alpha_j {z_j}}
\leq
(1+\delta_m)\sup_{m\leq j\leq n}\betr{\alpha_j}$
for all $m\leq n$.
If we have even 
$\sup_{j\leq n}(1-\delta_j)\betr{\alpha_j}
\leq \Norm{\sum_{j=1}^{n}\alpha_j {z_j}}
\leq \sup_{j\leq n}(1+\delta_j)\betr{\alpha_j}$
for all $n\in\N$ then $(z_j)$ is said to span $c_0$
asymptotically (or asymptotically isometrically).
While, by James' distortion theorem, a Banach space isomorphic to $c_0$ always
contains an almost isomorphic copy of $c_0$, Dowling, Johnson, Lennard and
Turett \cite{DJLT} proved the existence of a $c_0$-copy which does not contain
asymptotic copies of $c_0$.
Note that the L-structure of an L-embedded Banach space and,
respectively, the M-structure of its dual have an influence on the
existance of asymptotic copies of $\ell^1$ and, respectively, $c_0$.
For example,
it has been proved in \cite{Pfi-Fixp} that each almost isometric copy
of $\ell^1$ inside an L-embedded space contains an asymtotic copy of
$\ell^1$ (see \cite{Pfi-Fixp} also for the definitions)
and it has been proved there that if an L-embedded space is the dual
of an M-embedded space then its dual contains asymptotic copies of
$c_0$.\\
3.
There is an interesting difference of the construction of a $c_0$-copy
in the present proof and in the proof of property (V$^*$).
The latter one works for both separable and non-separable L-embedded
spaces whereas the present proof of property (X) runs into an obstacle
in the non-separable case:
Edgar \cite[Prop.\ 12]{Edg} showed that $\ell^1(\Gamma)$ has
property (X) if and only if $\mbox{card}(\Gamma)$ is not a real
measurable cardinality (that is if and only if there is no non-zero
measure on $\Gamma$ vanishing on singletons).
For a discussion of (X) and measurable cardinals
we refer to \cite{Neufang-X}.
It seems reasonable to conjecture
that an L-embedded Banach space may have property (X) if it does not
contain a subspace isomorphic to $\ell^1(\Gamma)$ with
$\mbox{card}(\Gamma)$ measurable, or, perhaps, if it has a dense subset
of non-measurable cardinality.\\
4.
Given a Banach space $Z$ it might occur that its bidual contains an
element $z_0^{**}$ which is L-direct to $Z$ that is
\bgl
\norm{z+z_0^{**}}=\norm{z}+\norm{z_0^{**}} \mbox{ for all } z\in Z.
\label{glst}
\egl
Godefroy has shown (\cite[IV.2]{HWW} or \cite[I.18.5.6]{JohLin})
that viewed as a function
on the unit ball of the dual $(B_{Z^*}, w^*)$ such an element $z_0^{**}$
is "very" discontinuous, for example it is nowhere continuous on
$(B_{Z^*}, w^*)$. The space $Z=\mbox{C}([0,1])$ and the function
$z_0^{**}=\charF{\Q\cap[0,1]}-\charF{(\R\setminus\Q)\cap[0,1]}$ serve
as an example. This function $z_0^{**}$ is of second Baire class but
does not belong to $Z$. In other words, the "local" property \Ref{glst}
is definitely weaker than the "global" one of being L-embedded because
if $Z$ were l-embedded then the second Baire class function $z_0^{**}$
would belong to $Z$ (cf.\ the second part of theorem \ref{theo-GT}).
\medskip\\
{\sc Acknowledgement}
I thank G. Godefroy and M. Neufang for interesting discussions,
and Dirk Werner for not only spotting but also correcting
a mistake in a previous version of the proof.
%\bibliography{Literatur}
%\bibliographystyle{plain}

\bigskip\noindent
Hermann Pfitzner\\
Universit\'e d'Orl\'eans\\
BP 6759\\
F-45067 Orl\'eans Cedex 2\\
France\\
e-mail: pfitzner@labomath.univ-orleans.fr
\end{document}